\documentclass[11pt]{amsart}
\usepackage{amssymb}
\usepackage{graphicx}
\usepackage{bigstrut}

\begin{document}
\newtheorem{theorem}{Theorem}[section]
\newtheorem{proposition}[theorem]{Proposition}
\newtheorem{lemma}[theorem]{Lemma}
\newtheorem{corollary}[theorem]{Corollary}
\newtheorem{defn}[theorem]{Definition}
\newtheorem{conjecture}[theorem]{Conjecture}

\theoremstyle{definition}
\newtheorem{definition}{Definition}

\numberwithin{equation}{section}

\newcommand{\R}{\mathbb R}
\newcommand{\TT}{\mathbb T}
\newcommand{\Z}{\mathbb Z}
\newcommand{\Q}{\mathbb Q}
\newcommand{\Sol}{\operatorname{Sol}}
\newcommand{\ND}{\operatorname{ND}}
\newcommand{\ord}{\operatorname{ord}}
\newcommand{\leg}[2]{\left( \frac{#1}{#2} \right)}  
\newcommand{\Sym}{\operatorname{Sym}}
\newcommand{\vE}{\mathcal E} 
\newcommand{\ave}[1]{\left\langle#1\right\rangle} 
\newcommand{\Var}{\operatorname{Var}}
 \newcommand{\tr}{\operatorname{tr}}
\newcommand{\supp}{\operatorname{Supp}}
\newcommand{\intinf}{\int_{-\infty}^\infty}

\newcommand{\dist}{\operatorname{dist}}
\newcommand{\area}{\operatorname{area}}
\newcommand{\vol}{\operatorname{vol}}
\newcommand{\diam}{\operatorname{diam}}

\newcommand{\Bin}{\operatorname{Binom}}
\newcommand{\pois}{\operatorname{Pois}}
\newcommand{\myK}{V}
\newcommand{ \myvar}{\operatorname{Var}}
\newcommand{\ripleyK}{\^K}
\newcommand{\E}{\mathbb{E}}
\newcommand{\Prob}{\operatorname{Prob}}

\renewcommand{\^}{\widehat}
\newcommand{\bx}{\mathbf x}
\newcommand{\by}{\mathbf y}
\newcommand{\bz}{\mathbf z}
\newcommand{\bl}{\mathbf \lambda}
\newcommand{\SF}{\mathcal H} 

\title{Local statistics of lattice points on the sphere}
\author{Jean Bourgain, Peter Sarnak and Ze\'ev Rudnick}
\date{\today}

\address{School of Mathematics, Institute for Advanced Study,
Princeton, NJ 08540 } \email{bourgain@ias.edu}

\address{Raymond and Beverly Sackler School of Mathematical Sciences,
Tel Aviv University, Tel Aviv 69978, Israel}
\email{rudnick@post.tau.ac.il}

\address{Department of Mathematics, Princeton University, Fine Hall,
Washington Road, Princeton, NJ 08544 and School of Mathematics,
Institute for Advanced Study, Einstein Drive, Princeton, NJ 08540
USA}

\begin{abstract}{A celebrated result of Legendre and Gauss determines
which integers can be represented as a sum of three squares, and for
those it is typically the case that there are many ways of doing so.
These different representations give collections of points on the
unit sphere, and a fundamental result, conjectured by Linnik, is
that under a simple condition these become uniformly distributed on
the sphere. In this note we survey some of our recent work, which
explores what happens beyond uniform distribution, giving evidence
to randomness on smaller scales. We treat the electrostatic energy,
local statistics such as the point pair statistic (Ripley's
function), nearest neighbour statistics, minimum spacing and
covering radius. We briefly discuss the situation in other
dimensions, which is very different. In an appendix we compute the
corresponding quantities for random points.}
\end{abstract}

\thanks{J.B. was supported in part by N.S.F. grants DMS-1301619
 and DMS 0835373.
Z.R. has received funding from the European Research Council under
the European Union's Seventh Framework Programme (FP7/2007-2013) /
ERC grant agreement n$^{\text{o}}$ 320755, and from the Israel
Science Foundation (grant No. 1083/10).
 P.S.  is partially supported by NSF
grants DMS-0758299 and DMS 1302952.
}


\maketitle

\section{Statement of results}

The set of integer solutions $(x_1,x_2,x_3)$ to the equation
\begin{equation}
x_1^2+x_2^2+x_3^2=n
\end{equation}
has been much studied. However it appears that the spatial
distribution of these solutions at small and critical scales as
$n\to\infty$ have not been addressed. The main results announced
below give strong evidence to the thesis that the solutions behave
randomly. This is in sharp contrast to what happens with sums of two
or four or more squares.

First we clarify what we mean by random. For a homogeneous space
like the $k$-dimensional sphere $S^k$ with its rotation-invariant
probability measure $\^\sigma$,  the binomial process is what you
get by placing $N$ points $P_1,\dots,P_N$ on $S^k$ independently
according to $\^\sigma$.
We are in interested in statistics, that is functions $f(P_1,\dots,
P_N)$, which have a given behaviour almost surely, as $N\to \infty$.
If this happens we say that this behaviour of $f$ is that of random
points. We shall also contrast features of random points sets with
those of ``rigid'' configurations, by which we mean points on a
planar lattice, such as the honeycomb lattice.

A celebrated result of Legendre/Gauss asserts that $n$ is a sum of
three squares if and only if $n\neq 4^a(8b+7)$. Let $ \vE(n)$ be the
set of solutions
\begin{equation}
 \vE(n) = \{\bx\in \Z^3: |\bx|^2=n\}
\end{equation}
and set
  \begin{equation}
  N= N_n:=\#\vE(n)\;.
  \end{equation}
The behaviour of $N_n$ is very subtle and it was a fine achievement
in the 1930's when it was shown that $N_n$ goes to infinity with $n$
(assuming say that $n$ is square-free; if $n=4^a$ then there are
only six solutions). It is known that $N_n\ll n^{1/2+o(1)}$ and if
there are {\em primitive} lattice points, that is
$\bx=(x_1,x_2,x_3)$ with $\gcd(x_1,x_2,x_3)=1$ (which happens if and
only if $n\neq 0,4,7\bmod 8$) then there is a lower bound of $
N_n\gg n^{1/2-o(1)}$. This lower bound is ineffective and indicates
that the behaviour of $N_n$ is still far from being understood
\cite{Si1}.

\begin{figure}[h]
\begin{center}
  \includegraphics[width=130mm]
{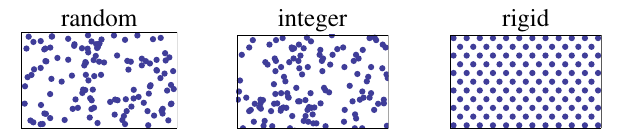}
 \caption{Lattice points coming from the
prime $n=1299709$ (center), versus random points (left) and rigid
points (right). The plot displays an area containing about $120$
points.} \label{random vs sphere pts}
\end{center}
\end{figure}

The starting point of our investigation is the fundamental result
conjectured by Linnik (and proved by him assuming the Generalized
Riemann Hypothesis), that for $n\neq 0,4,7 \bmod 8$, the points
\begin{equation}
 \^\vE(n):= \frac 1{\sqrt{n}}\vE(n) \subset S^2
\end{equation}
obtained by projecting to the unit sphere, become equidistributed on
the unit sphere with respect to $\^\sigma$  as $n\to \infty$. This
was proved unconditionally by Duke \cite{Duke, Duke-SP} and Golubeva
and Fomenko \cite{GF}, following a breakthrough by Iwaniec
\cite{Iwaniec}. Random points are equidistributed by definition and
the above result says that on this crudest global scale the
projected lattice points $\^\vE(n)$ behave like random points.
Figure~\ref{random vs sphere pts} gives some visual support for
random behaviour of  $\^\vE(n)$.

To make this precise we examine various statistics associated with
the placement of points in $S^2$. Our choice  of these statistics is
based on robustness tests for the random hypothesis, as well as
quantities which are of interest in number theoretical and harmonic
analysis applications. Our philosophy in what follows is that the
behaviour of a quantity in question is easy to determine for random
points (see  Appendix~\ref{sec:statistics}) while for $\^\vE(n)$ we
settle for estimates for them and also formulate conjectures,
 which are more precise. That one has to settle for such information
 for this kind of problem is to be expected given the problematic non-random
 behaviour of the number $N_n$ itself. The proofs of our assertions
 on $\^\vE(n)$ will appear in \cite{BRS}.

\subsection{Electrostatic energy}
The electrostatic energy of $N$ points $P_1,\dots,P_N$ on $S^2$ is
given by
\begin{equation}\label{electrostatic energy sum}
   E(P_1,\dots, P_N):=  \sum_{i\neq j} \frac 1{|P_i-P_j| }\;.
  \end{equation}
Here and in the sequel, $|x-y|$ is the Euclidean distance in $\R^3$.
This energy $E$ depends on both the global distribution of the
points as well as a moderate penalty for putting the points too
close to each other. The minimum energy configuration is known to
satisfy \cite{Wagner1, Wagner2, Brauchart}
\begin{equation}
 N^2-\beta N^{3/2}\leq \min_{P_1,\dots,P_N} E(P_1,\dots,P_N) \leq N^2-\alpha N^{3/2}
\end{equation}
for some $0<\alpha \leq \beta<\infty$. The configurations which
achieve this are rigid in various senses \cite{Dahlberg} and we will
see below in Corollary~\ref{cor:absolute continuity} that our points
$\^\vE(n)$ are far from being rigid. For random points one
has\footnote{Here and elsewhere, $\sim$ is the usual asymptotic
symbol denoting convergence to one of the ratio of the two sides.}
that $E\sim N^2$ but  the difference  $E-N(N-1)$ from the expected
value has no definite sign. Our first result is that to leading
order the points $\^\vE(n)$ have the same energy as random or minimal energy configurations.
\begin{theorem}\label{thm:electrostatic energy}
There is some $\delta>0$ so that
  \begin{equation}
   E(\^\vE(n)) =   N^2 +O(N^{2-\delta})
  \end{equation}
as $n\to \infty$, $n\neq 0,4,7 \mod 8$.
\end{theorem}
 We have not  been able to say anything about the sign of $E(\^\vE(n))-N(N-1)$ which according to
Table~\ref{table energies} appears to vary.

\begin{table}
\begin{tabular}{lrr}
\hline
&\multicolumn{2}{l}{$E-N(N-1)$} \\
\cline{2-3}
$N$ & integer & random   \bigstrut \\
\hline
$1224$  & $-282.$ & $95. $ \bigstrut  \\
 $3072$     & $37732.$     &  $-4704. $ \bigstrut   \\
$4296$   & $8380.$  & $ 1747. $ \bigstrut \\
\hline
\bigskip
\end{tabular}
\caption{The difference $E-N(N-1)$ between the electrostatic energy
  and its expected value, for various values of $N$. In the column
  labeled ``integer'', the energy for $\^\vE(n)$ was computed for the
  primes $n=104773$, $104761$ and $1299763$ with $N_n$ listed in the
  left-most column.   In the random case the result is a mean value
  of $20$ runs.}
  \label{table energies}
\end{table}

\subsection{Point pair statistics}
The point pair statistic and its variants are at the heart of our
investigation. It is a robust statistic as far as testing the
randomness hypothesis and it is called Ripley's function in the
statistics literature \cite{SKM}. For $P_1,\dots P_N\in S^2$ and
$0<r<2$, set
\begin{equation}\label{def of ripley}
 \ripleyK_r(P_1,\dots,P_N):=\sum_{\substack{i\neq j\\ |P_i-P_j|<r}} 1
\end{equation}
to be the number of ordered pairs of distinct points at (Euclidean)
distance at most $r$ apart. Note that we will allow $r$ to vary with $N$ so as to test randomness at different scales.

For fixed $\epsilon>0$, uniformly for
$N^{-1+\epsilon}\leq r\leq 2$, one has that for $N$ random points chosen with respect to the binomial process
\begin{equation}\label{int.9}
 \ripleyK_r(P_1,\dots,P_N) \sim \frac 14 N(N-1) r^2\;.
\end{equation}

Based on the results below as well as some numerical
experimentation, we conjecture that for $n$ square-free the points
$\^\vE(n)$ behave randomly w.r.t. Ripley's statistic at scales
$N_n^{-1+\epsilon}\leq r\leq 2$; that is
\begin{equation}\label{neweq1.10}
  \^K_r(\^\vE(n)) \sim \frac{N^2r^2}4, \quad \mbox{as } n\to \infty \;.
\end{equation}
 One of our main results is the
following which shows that \eqref{neweq1.10} is true at least in
terms of an upper bound which is off only by a multiplicative
constant.
\begin{theorem}\label{thm:poisson}
 Assume the Generalized Riemann Hypothesis (GRH). Then for fixed
 $\epsilon>0$ and $N^{-1+\epsilon} \leq r\leq 2$,
\begin{equation*}
 \ripleyK_r(\^\vE(n)) \ll_\epsilon N^2 r^2
\end{equation*}
where the implied constant depends only on $\epsilon$.
\end{theorem}
Remark: We do not need the full force of GRH here, but rather that
there are no ``Siegel zeros''.

We have not succeeded in giving individual lower bounds for
$\^K_r(\^\vE(n))$. What we can show is that at the smallest scale, that is $r$ of order $N_n^{-1+o(1)}\approx n^{-1/2+o(1)}$,
\eqref{neweq1.10} holds for most $n$'s:
\begin{theorem}\label{thm: lowerbd}
  There is some $\delta_0>0$ such that  for fixed
  $0<\delta<\delta_0$ and $r=n^{\delta-\frac 12}$,
  \begin{equation*}
     \^K_r(\^\vE(n)) \sim \frac{N^2r^2}4
  \end{equation*}
for almost all $n$.
\end{theorem}
  The constant $\delta_0$ can be determined explicitly, 
and is limited in our analysis by $h$ having to
be small in \eqref{eq 2.1} below . 

\subsection{Nearest neighbour statistics}
\begin{figure}[h]
\begin{center}
\includegraphics[width=100mm]
{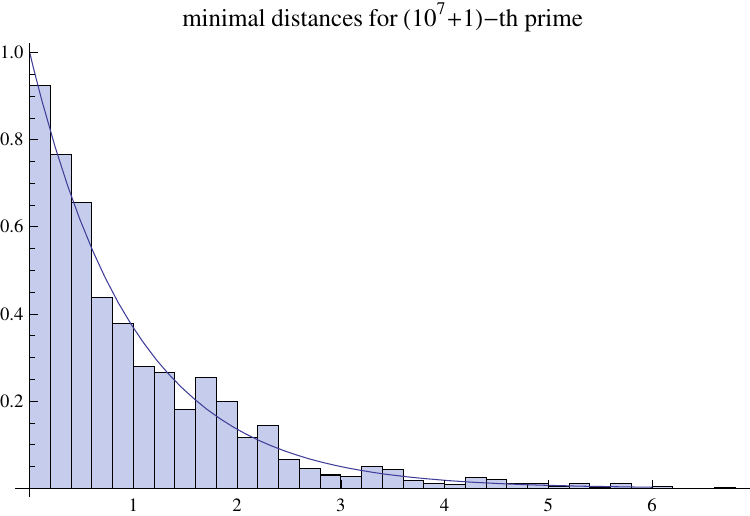}
\caption{A histogram of the scaled minimal spacing between lattice
points for for $n=179424691$, the $10,000,001$-th prime,  where
$N_n=94536$, and modulo symmetries there are $1970$ points.
The smooth curve is the exponential distribution $e^{-s}$.}
\label{hist10000001}
\end{center}
\end{figure}

Closely connected to $\ripleyK$ is the distribution of nearest
neighbour distances $d_j$, i.~e. the distance from $P_j$ to the
remaining points. It is more  convenient to work with these squares of the
distances.  Area considerations show that $\sum_j d_j^2 \leq
16$. For random points, the mean of $d_j^2$ is $4/N$.
In order to space these numbers at a scale for which they
have a limiting distribution in the random case, we rescale them by
their mean for the random case, i.e. replace $d_j^2$ by $\frac N4
d_j^2$. Thus for $P_1,\dots, P_N\in S^2$ define the nearest
neighbour spacing measure $\mu(P_1,\dots,P_N)$ on $[0,\infty)$ by
\begin{equation}
 \mu(P_1,\dots, P_N) := \frac 1N \sum_{j=1}^N \delta_{\frac N4 d_j^2}
\end{equation}
where $\delta_\xi$ is a delta mass at $\xi\in \R$. Note that the
mean of $\mu$ is at most $1$ and that for random points we have
\begin{equation}\label{int.17}
 \mu(P_1,\dots,P_N)\to e^{-x}dx,\quad \mbox{as } N\to \infty\;.
\end{equation}
Based on this and numerical experiments (see
figure~\ref{hist10000001}) we conjecture:
\begin{conjecture}
 As $n\to \infty$ along square-free integers, $n\neq 7\bmod 8$,
\begin{equation}
 \mu(\^\vE(n)) \to e^{-x}dx \;.
\end{equation}
\end{conjecture}

As a Corollary to Theorem~\ref{thm:poisson} we have
\begin{corollary}\label{cor:absolute continuity}
 Assume GRH. If $\nu$ is a weak limit of the $\mu(\^\vE(n))$ then $\nu$
 is absolutely continuous, in fact there is an absolute constant $c_4>0$ such that
\begin{equation}
 \nu\leq c_4dx\;.
\end{equation}
\end{corollary}

Corollary~\ref{cor:absolute continuity} implies that the
$\^\vE(n)$'s are not rigid for large $n$ since for rigid
configurations,
$\mu_{P_1,\dots,P_N} \to\delta_{\pi/\sqrt{12}}$. Also in as much as
it ensures that such a $\nu$ cannot charge $\{0\}$ positively, it
follows that almost all the points of $\^\vE(n)$ are essentially
separated with balls of radius approximately $N^{-1/2}$ from the
rest. Precisely, given a sequence $\eta_N$ satisfying
$\eta_N=o(N^{-1/2})$, all but $o(N)$ of the $N$ points in $\^\vE(n)$
have the ball of radius $\eta_N$ about them free of any other
points.

\subsection{Minimum spacing and covering radius}
Given $P_1,\dots,P_N\in S^2$ define the minimum spacing to be
\begin{equation}
m(P_1,\dots,P_N):=\min_{i\neq j} d_{i,j} = \min_j d_j \;.
\end{equation}
This statistic is very sensitive to the placement of points and it
is of arithmetic interest for $\^\vE(n)$. From the area packing
bound we have that
\begin{equation}
 m(P_1,\dots,P_N)\leq 4/\sqrt{N}
\end{equation}
for any configuration. In fact the rigid configuration of
Figure~\ref{random vs sphere pts} (coming from a planar lattice) maximizes $m$ asymptotically
\begin{equation}
  \max_{P_1,\dots,P_N} m(P_1,\dots,P_N) \sim \frac 2{\sqrt{N}}
  \cdot 2 \sqrt{ \frac{\pi}{\sqrt{12} }} \;.
\end{equation}
For random points the behaviour of the minimal spacing $m$ is very
different
\begin{equation}\label{int.23}
   m(P_1,\dots,P_N) = N^{-1+o(1)} \;.
\end{equation}
Based on the random point model as well as number theoretic
considerations which involve a nonlinear and shifted variation of
Vinogradov's least quadratic residue conjecture \cite{Vinogradov},
we pose
\begin{conjecture}\label{conj: minimal spacing}
  $m(\^\vE(n)) = N^{-1+o(1)}$ as $n\to \infty$.
\end{conjecture}
The lower bound in Conjecture~\ref{conj: minimal spacing} is an
immediate consequence of the integrality of the points in $\^\vE(n)$
since that implies that for the projected points $P_i\neq P_j \in
\^\vE(n)$, we have $|P_i-P_j|\geq 1/\sqrt{n}$ and since $N\geq
n^{1/2+o(1)}$ the lower bound follows. It is the upper bound that
appears difficult even assuming GRH.

As with the previous statistics we can establish the conjecture for
almost all $n$. Indeed it follows from Theorem~\ref{thm: lowerbd}
that
\begin{corollary}\label{corollary 1.8}
  Given $\epsilon>0$, $m(\^\vE(n)) \ll_\epsilon N^{-1+\epsilon}$ for
  almost all $n$.
\end{corollary}
Note that Conjecture~\ref{conj: minimal spacing} would follow from
the stronger conjecture of Linnik \cite{Linnikbook}, that for
$\epsilon>0$ and $n$ odd and square-free (and $n\neq 7\mod 8$) there
are $x_1,x_2,x_3$ with $|x_3|\leq n^\epsilon$ and
$x_1^2+x_2^2+x_3^2=n$, as such a representation provides a pair of
points $(x_1,x_2,\pm x_3)/\sqrt{n}\in \^\vE(n)$ at distance $\leq
n^{-1/2+o(1)}   \ll N^{-1+ \epsilon}$ from each other.

Finally we examine the covering radius for $\^\vE(n)$ though there
is little of substance that we can prove. Given $P_1,\dots , P_N\in
S^2$, the covering radius $M(P_1,\dots, P_N)$ is the least $r>0$ so
that every point of $S^2$ is within distance at most $r$ of some
$P_j$. Again an area covering argument shows that for any
configuration $M(P_1,\dots,P_N) \geq \frac 4{\sqrt{N}}$.

As a statistic, the covering radius $M$ is much more forgiving than
the minimal spacing $m$ in that the placement of a few bad points
does not affect $M$ drastically. In particular for random points,
$M\leq N^{-1/2+o(1)}$.

Based on this we conjecture the following, though admittedly with
much less evidence than the previous conjectures.
\begin{conjecture}
  $M(\^\vE(n)) = N^{-1/2+o(1)}$ as $n\to \infty$.
\end{conjecture}
An effective version of the equidistribution of $\^\vE(n)$ given in \cite{GF,
Duke-SP}, which is needed in the proof of
Theorem~\ref{thm:electrostatic energy}, yields an  $\alpha>0$ such
that $M(\^\vE(n)) \ll N^{-\alpha}$.

\subsection{Higher dimensions}

The distribution of the solutions to
\begin{equation}\label{sum of t squares}
x_1^2+x_2^2+\dots +x_t^2=n
\end{equation}
for $t\neq 3$ is very different and certainly non-random. Firstly
for $t=2$ and say $n$ a prime, $n=1\bmod 4$, there are exactly eight
solutions to \eqref{sum of t squares}. So there  is little to say
about the distribution for individual such $n$'s. However for
``generic" $n$'s which are sums of two squares, the projections of the
solutions  to the unit circle are uniformly distributed
\cite{Katai-Kornyei, Erdos-Hall}, and for such $n$'s the local
statistical questions certainly make sense.

For $t\geq4$, the projections onto the unit sphere of the solutions
to \eqref{sum of t squares} can be examined using the same
techniques that we use for $t=3$, with the main differences being
that the analysis is easier and the local behaviour is no longer
random. We only discuss the last feature and since it is only
enhanced with increasing $t$, we stick to $t=4$. Let $\vE_4(n)$ be
the set of solutions to \eqref{sum of t squares} and let
$\^\vE_4(n)$ be the projection of this set to $S^3$, the unit sphere
in $\R^4$. The first difference from $t=3$ is that
$N^{(4)}_n:=\#\^\vE_4(n)$ is a regularly behaved function of $n$.
When divided by $8$ it is multiplicative and for $n=p$ an odd prime
$\#\^\vE_4(p)=8(p+1)$. Thus the number of points $N=N_n^{(4)}$ being
placed on $S^3$ satisfies
\begin{equation}
  N=n^{1+o(1)}
\end{equation}
at least for odd $n$.
For $N$ random points on $S^3$ the two point function $\ripleyK_r$
defined as in \eqref{def of ripley} satisfies that for $\epsilon>0$
and $N^{-2/3+\epsilon}\leq r\leq 2$
\begin{equation}\label{int.25}
  \ripleyK_r(P_1,\dots,P_N) \sim N(N-1) V(r)
\end{equation}
where $V(r)$ is the relative volume of a cap $\{x\in S^3:
|x-x_0|<r\}$; for small $r$, $V(r)\sim \frac{2}{3\pi} r^3$.

On the other hand for $\^\vE_4(n)$, the integrality of the
corresponding points in $\vE_4(n)$, implies that for $x\neq y$,
$|x-y|\geq 1/\sqrt{n}$ and hence for $x\neq y\in \^\vE_4(n)$
\begin{equation}\label{Large spacings in dim 3}
  |x-y|\geq N^{-1/2+o(1)}\;.
\end{equation}
In particular for $r\leq N^{-1/2-\epsilon}$,
\begin{equation}
  \ripleyK_r(\^\vE_4(n)) = 0\;.
\end{equation}
Thus at the scales $N^{-2/3+\epsilon} \leq r\leq N^{-1/2-\epsilon}$,
the point pair function for $\^\vE_4(n)$ and that for random points
are very different.

This difference is also reflected in the minimum spacing function
$m(\^\vE_4(n))$  for $N$ points on $S^3$.
 From \eqref{Large spacings in dim 3} we have the
lower bound $m(\^\vE_4(n))\geq N^{-1/2+o(1)}$ and on the other hand
 there is a similar upper bound, namely
\begin{proposition}\label{prop:minspacing}
\begin{equation}
m(\^\vE_4(n)) = N^{-1/2+o(1)} \;.
\end{equation}
\end{proposition}
This is in sharp contrast to random points on $S^3$ for which
\begin{equation}\label{int.29}
  m(P_1,\dots,P_N) = N^{-2/3+o(1)}\;.
\end{equation}
Thus the points $\^\vE_4(n)$ are much more rigid than random points
but they are far from being fully rigid as the latter satisfy
(locally these points are placed at the vertices of the face
centered cubic lattice \cite{Conway-Sloane}):
\begin{equation}
  \max_{P_1,\dots,P_N} m_4(P_1,\dots,P_N) \sim \frac 2{N^{1/3}} c,
  \quad c=  \frac{\pi^{2/3}}{\sqrt{2}}\;.
\end{equation}

The nonrandom behaviour of the points $\^\vE_4(n)$ manifests itself
at a much larger scale as well, as is demonstrated by the minimum
covering radius $M_4(P_1,\dots,P_N)$. While being very nonrigid,
random points cover $S^3$ quite well. For them we have
\begin{equation}\label{int.31}
  M_4(P_1,\dots, P_N) = N^{-1/3+o(1)} \;.
\end{equation}
Somewhat surprisingly the points $\^\vE_4(n)$ which are more rigid
than random points, are poorly distributed in terms of covering.
This phenomenon of what might be called ``big holes'' was first
observed in the context of approximations of $2\times 2$ real
matrices by certain rational ones, by Harman \cite{Harman}. For
$\^\vE_4(n)$ we have
\begin{proposition}\label{prop:coveringradius}
\begin{equation}
  M(\^\vE_4(n)) \geq N^{-1/4+o(1)} \;.
\end{equation}
\end{proposition}

\section{Outline of the proofs}

For $n$ squarefree the general mass formula of Minkowski and Siegel
, which in the following special case is due to Gauss, expresses
$N_n$ in terms of $L(1,\chi_{d_n})$ where $\chi_{d_n}$ is the
quadratic character associated to the field $\Q(\sqrt{-n})$ of
discriminant $d_n$. From this and Siegel's lower bound on
$L(1,\chi_d)$ it follows that $N_n\gg n^{1/2-\epsilon}$ for any
$\epsilon>0$ (ineffectively). The key tool in our analysis of the
local point-pair functions is the mass formula applied to the
representations of the binary form $nu^2+2tuv +nv^2$ by the ternary
form $x_1^2+x_2^2+x_3^2 = \langle x, x\rangle$. Since this ternary
form has one class   in its genus, the mass formula gives the number
$A(n,t)$ of pairs $(x,y)\in \vE(n)\times \vE(n)$ with $\langle x,y
\rangle = t$, as  a product of local densities. Again this
is a special case of the mass formula, for which an elementary proof
as well as an explicit form was given in \cite{Venkov}, and this was
a critical ingredient in Linnik's approach to the equidistribution
of $\^\vE(n)$ (see \cite{EMV} for a recent exposition and extension
of his method). The local to global formula allows us to give rather
sharp upper bounds for $A(n,t)$. These are then used to control the
contributions of nearby points in the sum \eqref{electrostatic
energy sum} in the course of proving Theorem~\ref{thm:electrostatic
energy}. For pairs of points that are not too close we use modular
forms and in particular Duke's theorem. Specifically we effectivise
that analysis by giving a  power saving (namely $N^{-\alpha}$, for
some $\alpha>0$) upper bound for the spherical cap discrepancy of
the points $\^\vE(n)$. Putting these two together leads to
Theorem~\ref{thm:electrostatic energy}.

The proof of Theorem~\ref{thm:poisson} also uses the local formula
for $A(n,t)$, this time giving upper bounds for this quantity when
summed over $t$ in short intervals. It is critical that these upper
bounds are sharp up to a universal factor and depend only on the
subtle function $N_n$ and not on $n$. We achieve this by adapting
the upper bound sieve method of Nair \cite{Nair} to our setting.
This leads to an upper bound in terms of a product of local
densities of primes connected with $\chi_{d_n}$. It is here that we
need to assume that there are no Siegel zeros in order to ensure
that there is no dependence on $n$.

The ``almost all" result in Theorem~\ref{thm: lowerbd} is proven by
computing the asymptotic mean and variance of $\^K_r(\^\vE(n))
-N_n^2r^2/4$, with $n \leq R$. This is approached by analyzing
similar asymptotics for
\begin{equation}\label{eq 2.1}
K_h(\vE(n)) = \sum_{\substack{x,y\in \vE(n)\\x-y=h}}1
\end{equation}
and
\begin{equation}
K_{h,k}(\vE(n)) =\sum_{\substack{x,y,z,w\in \vE(n)\\x-y=h,\quad
z-w=k}}1 \;,
\end{equation}
($0\neq h,k\in \Z^3$).

The behaviour as $R\to \infty$ of $\sum_{n\leq R} K_h(\vE(n))$ may
be determined elementarily, while that of $\sum_{h\leq R}
K_{h,k}(\vE(n))$ can be derived using Kloosterman's circle method
for quadratic forms in $4$ variables (see for example \cite{Mal},
\cite{HB}).The leading terms are given as products of
Hardy-Littlewood local densities. The behaviour of $\sum_{n\leq R}
\^K_r(\vE(n)) N_n$ and $\sum_{n\leq R} N_n^2$ may be determined
using the Besicovich $r$-almost periodic properties of
$N_n/\sqrt{n}$ \cite{Peter}. We rederive this almost periodicity
directly using the circle method and this allows us to compare the
various local densities directly.

The proof of Proposition~\ref{prop:minspacing} is immediate from
Legendre and Gauss' Theorem. Namely $n-a^2=x_1^2+x_2^2+x_3^2$ has a
solution for $a=1$ or $a=2$ (recall $n$ is odd).
Proposition~\ref{prop:coveringradius} follows by considering annuli
about the north pole $(1,0,0,0)$.

\appendix

\section{Spatial statistics} \label{sec:statistics}
We give proofs of the statements that were made about the placement
of $N$ random points on $S^k$ and which were used to support the
thesis that our points $\vE(n)$  on $S^2$ behave like random points,
while in higher dimension the corresponding points are non-random.

We first prove  statements \eqref{int.9} and \eqref{int.25}
concerning Ripley's function. For $P_1,\dots,P_N\in S^k$,
\begin{equation}
\ripleyK_r(P_1,\dots,P_N) = \sum_{i\neq j} I_r(P_i,P_j)
\end{equation}
where $0\leq r\leq 2$ and $I_r(P_i,P_j)=1$ if $|P_i-P_j|\leq r$ and
is zero otherwise. For such $r$ let $B(P,r)$ denote the spherical
cap about $P$ consisting of all points $Q\in S^k$ such that
$|P-Q|\leq r$. Let $V(r)$ denote the $\^\sigma$-normalized surface
measure of $B(P,r)$, which is independent of $P$. In particular, in
dimension two, $V(r)= r^2/4$, and for $S^k$ in general $V(r)$ scales like $r^k$ for small $r$.

The claim  is that for $N^{-2/k+\epsilon}\leq r\leq 2$, as $N\to
\infty$
\begin{equation}\label{A.2}
\ripleyK_r(P_1,\dots,P_N) \sim N(N-1)V(r)
\end{equation}
in probability, by which we mean that  for each fixed $\epsilon>0$,
\begin{equation}
\Prob\left\{(P_1,\dots,P_N): \left| \frac{\ripleyK_r}{N(N-1)V(r)}-1
\right|>\epsilon\right \} \to 0
\end{equation}
as $N\to \infty$.

This follows in a standard way from Chebyshev's inequality once we
show that the expected value of $\ripleyK_r$ is
\begin{equation}\label{A.3a}
\E(\ripleyK_r) = N(N-1)V(r)
\end{equation}
and that its variance is
\begin{equation}\label{A.3b}
\Var(\ripleyK_r) = 2N(N-1)V(r)(1-V(r))\;.
\end{equation}
Indeed, then
\begin{equation}
\E\left( (\frac{\ripleyK_r}{\E(\ripleyK_r)}-1 )^2 \right) =
\frac{2(1-V(r))}{N(N-1)V(r)}
\end{equation}
which tends to zero if and only if $N(N-1)V(r)\to \infty$. Our lower
bound  for $r$ ensures that the latter holds and if we allow $r$ to
be $N^{-2/k}$ or smaller then it is clear that \eqref{A.2} is no
longer valid.

To see that \eqref{A.3a} and \eqref{A.3b} hold note that
\begin{equation}
\int_{S^k}\int_{S^k} I_r(P_1,P_2) d\^\sigma(P_1)d\^\sigma(P_2) =
V(r)\;,
\end{equation}
while
\begin{equation}
\int_{S^k}\int_{S^k} \int_{S^k} I_r(P_1,P_2)I_r(P_2,P_3)
d\^\sigma(P_1)d\^\sigma(P_2) )d\^\sigma(P_3) = V(r)^2
\end{equation}
and
\begin{equation}
\int_{S^k}\int_{S^k}\int_{S^k} \int_{S^k} I_r(P_1,P_2)I_r(P_3,P_4)
d\^\sigma(P_1)d\^\sigma(P_2) )d\^\sigma(P_3)d\^\sigma(P_4) = V(r)^2\;.
\end{equation}
Hence
\begin{equation}
\E(\ripleyK_r) = \sum_{i\neq j} \E\Big(I_r(P_i,P_j)\Big) = N(N-1) V(r)\;,
\end{equation}
which gives \eqref{A.3a}, and
\begin{equation}
\begin{split}
\E(\ripleyK_r^2) & = \E\Big((2\sum_{i<j} I_r(P_i,P_j))^2\Big) \\
&=4\sum_{i_1<j_1}\sum_{i_2<j_2} \E\Big(I_r(P_{i_1},P_{j_1}) I_r(P_{i_2},P_{j_2}) \Big) \\
&=4 \sum_{\substack{ i_1=i_2\\j_1=j_2\\i_1<j_1}} V(r) + 4 \sum_{\substack{ (i_1,j_1)\neq (i_2,j_2)\\i_1<j_1\\i_2<j_2}} V(r)^2\\
&= 4V(r)\frac{N(N-1)}2 +4V(r)^2 \left( (\frac{N(N-1)}2)^2-\frac{N(N-1)}2 \right)\\
&= V(r)^2 \Big(N(N-1)\Big)^2 +2N(N-1)V(r)\Big(1-V(r)\Big)\;,
\end{split}
\end{equation}
which immediately gives \eqref{A.3b}.

Next, focusing on the two-dimensional case of $S^2$, we compute the expected value of the average $\frac
1N\sum_{j=1}^N d_j^2$
and show that it equals
\begin{equation}\label{mean d2}
\E(\frac 1N\sum_{j=1}^N d_j^2)=\frac 4N \;.
\end{equation}
  Note that $d_1,\dots, d_N$ are not {\em independent}, for instance for $N=2$ we clearly have $d_1=d_2$.
However, they do have the same distribution and hence
\begin{equation}\label{ex mean in terms of 1}
 \E(\frac 1N\sum_{j=1}^N d_j^2) = \E(d_1^2) =
\E\Big(\min_{j=2,\dots,N}|P_1-P_j|^2 \Big)\;.
\end{equation}
Now $|P_j-P_1|^2$   (for $j=2,\dots,N$) are i.i.d. and take values
in $[0,4]$, and for $0\leq x\leq 4$,
\begin{equation}
\begin{split}
\Prob\Big(\min_{j=2,\dots,N}|P_1-P_j|^2>x\Big) &=
\Prob\Big(|P_2-P_1|^2>x\Big)^{N-1}\\
&= \Big(1-V(\sqrt{x})\Big)^{N-1} = (1-\frac x4)^{N-1}\;.
\end{split}
\end{equation}
From general principles, if $Y$ is non-negative then
 $ \E(Y)  = \int_0^\infty \Prob(Y>y)dy$. Hence
\begin{equation}
\E\Big(\min_{j=2,\dots,N}|P_1-P_j|^2 \Big)  = \int_0^4(1-\frac x4)^{N-1}dx=
\frac 4N
\end{equation}
which in conjunction with \eqref{ex mean in terms of 1} proves
\eqref{mean d2}.

 We turn to \eqref{int.17}, the distribution of scaled nearest neighbour spacings 
\begin{equation}
\mu_N(P_1,\dots,P_N) = \frac 1N \sum_{j=1}^N
\delta_{\frac{Nd_j^2}4}\;.
\end{equation}
For $x\geq 0$ we examine the expectations
\begin{equation}
\E\Big(\mu_N[0,x]\Big) = \E \Big(\frac 1N \sum_{j=1}^N I(\frac{Nd_j^2}4 \leq
x)\Big)\;,
\end{equation}
where $I(\bullet)=1$ if the condition $\bullet$ holds, and $0$ otherwise.
Setting $r=2\sqrt{\frac xN}$ we have
\begin{equation}
\mu_N[0,x] = \frac 1N \sum_j I\Big(\min_{k\neq j} d(P_k,P_j) \leq r\Big)\;.
\end{equation}
Hence
\begin{equation}\label{A.11}
\mu_N[0,x] \leq A_1:= \frac 1N \sum_j \sum_{k\neq j} I_r(P_j,P_k)
\end{equation}
and
\begin{equation}
\mu_N[0,x] \geq A_1-A_2\;,
\end{equation}
where
\begin{equation}\label{A.12}
A_2 := \frac 1N \sum_j \sum_{\substack{ k_1,k_2\\ k_1\neq
k_2\\k_1,k_2\neq j}} I_r(P_j,P_{k_1}) I_r(P_j,P_{k_2})
\end{equation}
and continuing with inclusion/exclusion, defining $A_k$ analogously,
\begin{equation}\label{A.13}
\begin{split}
\mu_N[0,x] &\leq A_1-A_2 +A_3+\dots+A_{2\ell+1} \\
\mu_N[0,x] &\geq A_1-A_2 +A_3+\dots-A_{2\ell+2} \;.
\end{split}
\end{equation}
Taking expectations we have from \eqref{A.11} and \eqref{A.12}
\begin{equation}
\E\Big(\mu_N[0,x]\Big) \leq \frac 1 N N(N-1) V(r) \to x, \quad \mbox{ as }
N\to \infty
\end{equation}
and hence
\begin{equation}
\limsup_{N\to \infty}\E\Big( \mu_N[0,x]\Big) \leq x
\end{equation}
while
\begin{equation}
\E\Big(\mu_N[0,x]\Big) \geq \frac{N(N-1) V(r)}{N} - \frac 1N
\frac{(N-1)(N-2)}2 V(r)^2 \to x-\frac {x^2}2
\end{equation}
and hence
\begin{equation}
\liminf_{N\to \infty}\E\Big(\mu_N[0,x]\Big) \geq x-\frac {x^2}2 \;.
\end{equation}
Continuing with \eqref{A.13} taking expectations and limit as $N\to
\infty$ yields after a similar calculation that
\begin{equation}
\lim_{N\to \infty}  \E\Big(\mu_N[0,x]\Big) = x-\frac{x^2}{2}+\frac{x^3}{3!}
-\dots = 1-e^{-x}
\end{equation}
that is
\begin{equation}
\lim_{N\to \infty} \E(\mu_N) = e^{-t}dt \;.
\end{equation}

One can compute variances as we did in \eqref{A.3b} above, from
which it follows that
\begin{equation}
\mu_N[0,x] \to 1-e^{-x}
\end{equation}
in probability, which is what is meant in \eqref{int.17}.

We end this short appendix with proofs of \eqref{int.23},
\eqref{int.29} and \eqref{int.31}. The first two are concerned with
the minimum spacing for random points in $S^k$. The probability of
placing $N$ independent points in $S^k$ so that none are closer to
each other than $r$ is clearly
\begin{equation}\label{A.18}
\Big(1-V(r)\Big)\Big(1-2V(r)\Big)\dots \Big(1-(N-1)V(r)\Big)
\end{equation}
as long as $(N-1)V(r)<1$ (otherwise the probability is zero). From
this it is clear that if $r\geq N^{-2/k+\epsilon}$ so that $V(r)N^2
\geq N^{\epsilon'}$, then the product in \eqref{A.18} tends to $0$.
So that with probability tending to $1$ the minimum spacing is at
most $N^{-2/k+\epsilon}$. On the other hand of $r\leq
N^{-2/k-\epsilon}$ then $V(r)N^2 \leq N^{-\epsilon'}$ and the
product in \eqref{A.18} goes to $1$, so that with probability
tending to $1$, the minimal spacing is at least $N^{-2/k-\epsilon}$.
That is $m(P_1,\dots,P_N) =N^{-2/k+o(1)}$ with probability tending
to $1$, which establishes \eqref{int.23} and \eqref{int.29}.

Finally we turn to \eqref{int.31}, which is concerned with the
covering radius $M(P_1,\dots ,P_N)$ for random points. For any
configuration of points on $S^k$, surface area considerations show
that
\begin{equation}
M(P_1,\dots,P_N) \gg N^{-1/k}\;.
\end{equation}
So to establish \eqref{int.31} we need to show that for $\epsilon>0$
and random points $P_1,\dots,P_N$,
\begin{equation}
M(P_1,\dots,P_N) \ll N^{-1/k+\epsilon}\;.
\end{equation}
Given $x\in S^k$, the probability that the cap $B(x,r)$ does not
contain any of the points $P_1,\dots, P_N$ is
\begin{equation}
\Big(1-V(r)\Big)^N\;.
\end{equation}
Hence if $x_1,\dots, x_L$ are $L$ points in $S^k$, the probability
that at least one of the caps $B(x_j,r)$  does not contain any of
$P_1,\dots,P_N$ is at most
\begin{equation}\label{A.21}
L\Big(1-V(r)\Big)^N\;.
\end{equation}
In particular if $L=N$   and $r=N^{-1/k+\epsilon}$ then the
probability in \eqref{A.21} goes to $0$ as $N\to \infty$. Hence with
probability tending to $1$ each of the $N$ caps $B(x_j,r)$ contains
at least one of the $P_j$'s. Now choose the $x_j$'s to be a $c_k
N^{-1/k}$ cover (i.e. each point of $S^k$ is within $c_k N^{-1/k}$
of one of
 $x_1,\dots ,x_N $).
Since there is a $P_i$ in each cap $B(x_j,r)$, it follows that the
$P_i$'s form a $2N^{-1/k+\epsilon}$ covering of $S^k$. This
completes the proof of \eqref{int.31}.

\end{document}